\documentclass[11pt]{amsart}

\usepackage[dvipsnames]{xcolor}
\usepackage{amssymb,verbatim}
\usepackage{amsmath,amsfonts,enumitem,bm}
\usepackage[mathscr]{euscript} 
\usepackage{amsthm}
\usepackage{url}
\usepackage{graphicx} 
\usepackage{enumitem} 
\usepackage{hyperref} 
\hypersetup{colorlinks} 
\usepackage{bbm} 
\usepackage{bm} 
\usepackage{mathrsfs} 
\usepackage{cancel} 
\usepackage{ytableau} 
\usepackage{mathtools}
\usepackage{faktor} 
\usepackage{algorithm}
\usepackage{algorithmic}

\usepackage[colorinlistoftodos]{todonotes}

\RequirePackage{cleveref}
\usepackage{hypcap}
\hypersetup{colorlinks=true, citecolor=darkblue, linkcolor=darkblue}
\definecolor{darkblue}{rgb}{0.0,0,0.7}
\newcommand{\darkblue}{\color{darkblue}}

\definecolor{darkred}{rgb}{0.68,0,0}

\definecolor{darkgreen}{rgb}{0,.38,0}

\newcommand{\defn}[1]{\emph{\darkblue #1}}

\setlist[enumerate]{
	label=\textnormal{({\roman*})},
	ref={\roman*}}

\makeatletter
\def\th@plain{%
	\thm@notefont{}
	\itshape 
}
\def\th@definition{%
	\thm@notefont{}
	\normalfont 
}
\makeatother

\newtheorem{thm}{Theorem}[section]
\newtheorem{lemma}[thm]{Lemma}

\newtheorem*{claim*}{Claim}
\newtheorem{cor}[thm]{Corollary}

\newtheorem{prop}[thm]{Proposition}
\newtheorem{conj}[thm]{Conjecture}
\newtheorem{conv}[thm]{Convention}
\newtheorem{question}[thm]{Question}

\theoremstyle{definition}
\newtheorem{ex}[thm]{Example}

\numberwithin{figure}{section}
\numberwithin{equation}{section}

\def\emp{\nothing}

\def\nn{\mathbb N}
\def\cc{\mathbb C}

\def\kk{\mathbb K}

\def\la{\lambda}

\def\al{\alpha}

\def\<{\langle}
\def\>{\rangle}

\def\GL{ {\text {\rm GL} } }

\def\0{{\mathbf 0}}

\def\nothing{\varnothing}

\def\.{\hskip.06cm}
\def\ts{\hskip.03cm}

\def\.{\hskip.06cm}
\def\ts{\hskip.03cm}

\def\nin{\noindent}

\newcommand{\textsu}[1]{\textup{\textsf{#1}}}
\usepackage{indentfirst}
\DeclareTextSymbolDefault{\ae}{T1}

\newcommand{\ComCla}[1]{\textup{\textsu{#1}}}

\newcommand{\NP}{\ComCla{NP}}

\newcommand{\BPP}{\ComCla{BPP}}

\newcommand{\coRP}{\ComCla{coRP}}

\renewcommand{\P}{\ComCla{P}}

\def\poly{{\P}}

\newcommand{\inv}{\operatorname{{\ell}}}
\newcommand{\Des}{\operatorname{{\rm Des}}}

\newcommand{\Sch}{\mathfrak{S}} 

\newcommand{\Sc}{\mathfrak{S}}

\newcommand{\code}{{\sf{code}}}

\newcommand{\satA}[1]{*_{#1}}

\newcommand{\setsize}{\alpha}

\begin{document}

\title[Saturation property fails for Schubert coefficients]{Saturation property fails for Schubert coefficients}


\author[Igor Pak \. \and \. Colleen Robichaux]{Igor Pak$^\star$  \. \and \.  Colleen Robichaux$^\diamond$}

\thanks{\today}
\thanks{\thinspace ${\hspace{-.45ex}}^\star$Department of Mathematics,
UCLA, Los Angeles, CA 90095, USA. Email:  \texttt{pak@math.ucla.edu}}
\thanks{\thinspace ${\hspace{-.45ex}}^\diamond$Department of Mathematics,
UC Davis, Davis, CA 95616, USA. Email:  \texttt{robichaux@ucdavis.edu}}

\begin{abstract}
The saturation property for Littlewood--Richardson coefficients
was established by Knutson and Tao in 1999.  In 2004, Kirillov
conjectured that the saturation property extends to Schubert
coefficients.  We disprove this conjecture in a strong form,
by showing that it fails for a large family of instances.
\end{abstract}

\maketitle

\section{Introduction}\label{s:intro}

\subsection{Saturation property}\label{ss:intro-sat}
The \emph{saturation conjecture} (now \emph{saturation theorem})
was proven by Knutson and Tao \cite{KT99} and
is one of the most celebrated results in Algebraic Combinatorics.
Despite its apparent simplicity, it became the last piece of the puzzle resolving \emph{Horn's problem}, which describes possible spectra of
three Hermitian matrices satisfying the equation $A+B=C$.  See e.g.\
\cite{Buch00,Ful00} for overviews of different aspects
of this remarkable story and \cite{BVW17,Kum14} for some later developments.

The \defn{saturation theorem} \ts states that for all partitions \ts $\la,\mu,\nu$ \ts
with \ts $|\la|=|\mu|+|\nu|$, we have:
\begin{equation}\label{eq:saturation}
c^{\la}_{\mu\nu} \. > \. 0 \quad \Longleftrightarrow \quad c^{N\la}_{N\mu, N \nu} \. > \. 0 \ \ \text{\em for any \ $N\ge 1 \ts$},
\end{equation}
where \ts $c^\la_{\mu\nu}$ \ts denote the \emph{Littlewood--Richardson {\rm (LR)} coefficients},
see e.g.\ \cite[$\S$1.9]{Mac95}.  We refer to \eqref{eq:saturation} as the \emph{saturation property}.

Note that \. $\Rightarrow$ \. is the easy direction, a special case of the 
\emph{semigroup property}, that triples of partitions \ts $(\la,\mu,\nu)$ \ts
positive LR coefficients form a semigroup, see e.g.~\cite{Zel99}.  
On the other hand, the direction  \. $\Leftarrow$ \. says that this semigroup 
is saturated.  This is a difficult result which does not immediately follow from 
the \emph{Gelfand--Tsetlin pattern} combinatorial interpretation of LR~coefficients, 
the way that \. $\Rightarrow$ \. does.  
Known proofs involve technical combinatorial \cite{KT99}, algebraic \cite{DW00,KM08},
or algebro-geometric arguments \cite{Bel06}.

The LR coefficients play a central role in Algebraic Combinatorics and
related areas, so in the aftermath of the saturation theorem, a number
of generalizations of \eqref{eq:saturation}
have been proposed, see e.g.\ a large collection in~\cite{Kir04}.
Unfortunately, in the quarter century since the original paper, very
few saturation type properties have been established, all of them remarkable.
These include the quantum version by Belkale \cite[$\S$4.1]{Bel08}, the
general reductive group version by
Kapovich--Millson \cite{KM08} (see also \cite{BK10}), and
the equivariant version by Anderson--Richmond--Yong \cite{ARY13}
(see also~\cite{ARY19}).
Most recently, an unexpected M\"obius strip version by Min~\cite{Min24}
resolved the \emph{Gao--Orelowitz--Yong saturation conjecture} \ts
for the \emph{Newell--Littlewood numbers} \cite{GOY21}.

Among positive results, let us also mention
\emph{Fulton's conjecture} \ts resolved by Knutson, Tao, and Woodward \cite{KTW04},
which can be viewed as a variation on \eqref{eq:saturation}~$:$
\begin{equation}\label{eq:Fulton}
c^{\la}_{\mu\nu} \. = \. 1 \quad \Longleftrightarrow \quad c^{N\la}_{N\mu, N \nu}
\. = \. 1  \ \ \text{\em for any \ $N\ge 1 \ts$}.
\end{equation}
Here the direction \. $\Leftarrow$ \. follows easily from the saturation property,
while the direction \. $\Rightarrow$ \. requires further work.  This
\emph{uniqueness property} \ts was also generalized in several ways, notably in
\cite{BKR12}.

In the negative direction, there is a large number of counterexamples
to the saturation property for various extensions of LR coefficients that are
scattered across the literature.  For example, it was shown by \`Elashvili \cite{Ela92}
that the saturation fails for root system~$B$ (explaining the factor of~$2$ in \cite{BK10,KM08}), see also \cite[p.~340]{Zel99} and \cite[$\S$4.2]{DM06}.
Similarly, Buch observed that the saturation property easily fails for the $K$-theoretic generalization
\cite[p.~71]{Buch02}.  Most recently, the second author, Yadav, and Yong noted that the saturation
fails for Schur $P$-polynomials \cite[Remark~7.7]{RYY22}; see also \cite[$\S$7.4]{CR23} for a
larger example in the geometric context.

For \emph{Kronecker coefficients}, which famously generalize the LR~coefficients, 
the semigroup property was established in \cite{CHM07}.  However, the saturation 
property fails already for the two row partitions, see e.g.\ \cite[Ex.~2]{Ros01}.
Because of their crucial role in \emph{Geometric Complexity Theory} (GCT),
Mulmuley stated a weak version of the property \cite[$\S$1.6]{Mul09}.
Soon after, Briand, Orellana and Rosas  \cite{BOR09},
disproved Mulmuley's weak version.  Curiously, a version of
the uniqueness property \eqref{eq:Fulton} continues to hold for
Kronecker coefficients, see \cite{SS16} and references therein.

Finally, for \emph{reduced Kronecker coefficients}, the saturation property
was conjectured independently by Kirillov \cite[Conj.~2.33]{Kir04}
and Klyachko \cite[Conj.~6.2.4]{Kly04}. These constants occupy
an intermediate place between LR and Kronecker coefficients, as they
generalize the former and are a special case of the latter.
Only recently, the first author and Panova constructed
a large family of counterexamples in this case \cite{PP20};
see also~$\S$\ref{ss:finrem-sat-other} for the
computational complexity context.

\smallskip

\subsection{Schubert saturation}\label{ss:intro-schub}
\defn{Schubert coefficients} \ts 
$\{c_{u,v}^w \. : \. u,v,w\in S_\infty\}$ \ts
can be defined as structure constants for \emph{Schubert polynomials}:
\begin{equation}\label{eq:Schub-def}
\Sc_u\. \cdot\.\Sc_v\ = \ \sum_{w\ts\in\ts S_\infty}\. c_{u,v}^w\. \Sc_w\..
\end{equation}
Here $S_\infty$ consists of permutations with all but finitely many fixed points,
and \ts $\Sc_u \in \nn[x_1,x_2,\ldots]$.

It is known that \ts $c_{u,v}^w$ \ts are nonnegative integers, as they
count certain intersection numbers.
These coefficients play a major role in \emph{Schubert calculus}, a rapidly developing
area of algebraic geometry, motivated in part by rich connections with
representation theory and algebraic combinatorics, see e.g.\ \cite{AF24,Knu22}.

We note that \emph{Grassmannian permutations}~$w$ (permutations with at most one descent)
correspond to integer partitions \ts $\la=\la(w)$. Famously, the corresponding Schubert
polynomials are symmetric and coincide with Schur polynomials: $\Sch_w=s_{\la(w)}$.
Thus, Schubert coefficients can be viewed as advanced generalizations of LR~coefficients, see e.g.\ \cite[$\S$2.6.4]{Man01}.

It is then natural to ask  if Schubert coefficients also satisfy a saturation property
extending \eqref{eq:saturation}.  In~\cite{Kir04}, Kirillov formulated
this as a conjecture (see below), which remained open until now.  Motivated
by GCT and with complexity applications as a motivation
(cf.~$\S$\ref{ss:finrem-CS}), Mulmuley also speculated that saturation might hold in this case \cite[$\S$3.7]{Mul09}.

For a permutation $w\in S_n\ts$, the \defn{Lehmer code}, also called
the \emph{inversion index}, is defined as
$$
\code(w) \, := \, (c_1,\ldots,c_n)\in \nn^n, \quad \text{where} \quad c_i \. := \.
|\{j\,: \, j>i, \. w(j)<i\}|.
$$
Clearly, \ts $c_1+\ldots+c_n = \inv(w)$ \ts is the \emph{number of inversions}.
We can now define the operation of \defn{code scaling} \ts as follows:
\[N\satA{} w \, := \, \code^{-1}(N c_1,N c_2,\ldots,N c_n,0,\ldots,0) \ts \in \ts S_{Nn}\..
\]

It is  easy to see that the code of a Grassmannian
permutation $w$ is a partition $\la(w)$ written in reverse,
so this code scaling corresponds to the usual multiplication
of partitions by a constant~$N$.
Therefore, the following conjecture is a natural generalization
of the saturation theorem:

\smallskip

\begin{conj}[{\rm  Kirillov \cite[Conj.~6.28]{Kir04}}{}] \label{conj:Kir}
For every \ts $u,v,w\in S_n\ts$, we have:
\begin{equation}\label{eq:sat-Kir}
c^{w}_{u,v} \. > \. 0 \quad \Longleftrightarrow \quad c_{N\ast u,\ts N\ast v}^{N\ast w}
\. > \. 0  \ \ \text{\em for any \ $N\ge 1 \ts$}.
\end{equation}
\end{conj}

\smallskip

We disprove Kirillov's conjecture in a strong form, by constructing
a large family of triples of permutations for which the saturation
property fails:

\smallskip

\begin{thm}\label{t:sat-Kir}
Let \. $u\in S_n$ \. and let \. $j-i\geq 2$, such that:
\begin{enumerate}
    \item $u\. \lessdot \. ut_{ij}$,
    \item the set $\{c \. : \. u(i)<c<u(j), \. j<u^{-1}(c) \}\neq\emp$, and
    \item if $i<k<j$, then $u(k)<u(j)$,
\end{enumerate}
where \. $t_{ij}=(i,j)$ \. is a transposition and \. $x\lessdot y$ \.
is the cover relation in the strong Bruhat order of~$S_n\ts$.
Finally, let \. $v=(i,i+1)$ \. be a simple transposition, and let \. $w=ut_{ij}$. Then:
\begin{equation}\label{eq:sat-Kir-thm}
c_{u,v}^w\.= \. 1 \qquad \text{and}  \qquad c_{N\satA{}u,N\satA{}v}^{N\satA{}w}\. = \. 0 \ \  \text{for all}  \ \ N>1\ts.
\end{equation}
\end{thm}

\smallskip

We illustrate the theorem
by an explicit sequence of counterexamples:

\smallskip

\begin{cor}
\label{c:main-Kir}
For all \ts $n\ge 4$, let \ts $v:=(n-2,n-3)$ \ts be a simple transposition in
\ts $S_{n}\ts$, and let
$$
u := (2,3,\ldots,n-2,1,n,n-1)\., \qquad
w := (2,3,\ldots,n,1,n-2,n-1)\ts.
$$
Then:
$$
c_{u,v}^w \. = \. 1 \qquad \text{and} \qquad c_{N\ast u,\ts N\ast v}^{N\ast w} \. = \. 0 \ \ \ \text{for all} \ \ N>1\ts.
$$
\end{cor}

\smallskip

In particular, this shows that Kirillov's Conjecture~\ref{conj:Kir} fails
for all \ts $n\ge 4$ \ts and all \ts $N \ge 2$.
Note also that in Theorem~\ref{t:sat-Kir}, we always have \ts $v$ \ts is
a simple transposition.  In this case,
Schubert coefficients have a simple combinatorial interpretation
given by \emph{Monk's rule} (Proposition~\ref{prop:monk}).  The proof of
Theorem~\ref{t:sat-Kir} (see~$\S$\ref{ss:Kir-proof}) uses this
combinatorial rule for the first part of \eqref{eq:sat-Kir-thm}, and the
\emph{St.~Dizier--Yong vanishing condition}
(Lemma~\ref{thm:stDY}) for the second part.

Let us emphasize that although \eqref{eq:sat-Kir}
is a direct generalization of the saturation property \eqref{eq:saturation},
it fails for what was originally an ``easy direction'' \. $\Rightarrow\.$, i.e.\ 
triples of permutations \ts $(u,v,w)$ \ts encoded by their Lehmer's codes 
do not form a semigroup.\footnote{After this paper appeared, 
the first author and Slonim refuted the ``hard direction'' \. $\Leftarrow\.$ as well, see \cite[Prop.~1.2]{PS26}.}
Also, observe that
the permutations in Corollary~\ref{c:main-Kir} have at most $2$ descents.
In this case, Schubert coefficients have two different combinatorial interpretations
\cite{Cos09,BKPT16},
cf.~$\S$\ref{ss:finrem-CS}.

\smallskip



\subsection{Structure of the paper} \label{ss:intro-structure}
We start with the algebraic combinatorics background in Section~\ref{s:back},
where we include standard definitions, notation and basic results in the area.
In Section~\ref{s:Kir}, we prove Theorem~\ref{t:sat-Kir}, thus giving counterexamples
to Kirillov's Conjecture~\ref{conj:Kir}.  We conclude with final remarks and
open problems in Section~\ref{s:finrem}.  Notably,  we discuss the failure of
the saturation property \eqref{eq:sat-Kir} in the context of the complexity of
Schubert vanishing.

\medskip



\section{Background}\label{s:back}

\subsection{Basic notation}\label{ss:back-basic}
We use \ts $\nn=\{0,1,2,\ldots\}$ \ts and \ts $[n]=\{1,2,\ldots,n\}$.
To simplify the notation, for a set $A$ and elements $x,y$,
we write \ts $A+x$ \ts to denote \ts $A\cup\{x\}$,
and \ts $A-y$ \ts to denote \ts $A-\{y\}$.

We use \ts $\kk$ \ts  to denote the set of infinite sequences with entries in $\nn$ with
finite support.  When writing such sequences, we omit the infinite
tail of trailing zeros, and write only the prefix with the support of the sequence,
so e.g.\ \ts $30120000\ldots$ \ts is written as \ts $3012$.

We use $S_n$ to denote the symmetric group, which we view as the group of
permutations of~$[n]$.  Denote by \ts $\iota$ \ts
the inclusion $\iota: S_n \. \hookrightarrow \. S_{n+1}$ \ts defined by
$w(1)\cdots \. w(n) \, \mapsto \, w(1) \ \cdots \. w(n) \,\. n+1$.
As above, let \ts $S_\infty=\bigcup_{n\geq 1} S_n$ \ts
denote the group of permutations of \ts $\nn_{\ge 1}=\{1,2,\ldots\}$ \ts with all but finitely
many fixed points, where the inclusion is given by~$\iota$.
By analogy with infinite sequences, when writing  permutations
we omit the tail of fixed points, so e.g.\ $34125678\ldots$ \ts is written as $3412$.

For convenience, we use both the sequence and word notation for permutations,
so for example \ts $(3,1,4,2)$ \ts and \ts $3142$ \ts correspond to the same
permutation in~$S_4\ts$.  We also use \ts $s_i=(i,i+1)$ \ts
to denote simple transpositions swapping $i$ and $i+1$,
and \ts $t_{ij}=(i,j)$ \ts to denote  general transpositions
swapping $i$ and $j$, where \ts $i<j$.  We hope this does not lead to confusion.

For a permutation $w\in S_n\ts$, let \ts $\ell(w):=|\{(i,j)\.:\. 1\le i<j\le n, \. w(i)>w(j)\}|$
\ts denote the \emph{number of inversions} \ts in~$w$.  For permutations $u,v \in S_n\ts$,
we write \ts $u \lessdot v$ \ts if \ts $v = u t_{ij}$ \ts for some $i<j$, and $\ell(v)=\ell(u)+1$.
This is the \emph{cover relation} \ts for the \emph{strong Bruhat order}, which can now be defined
by transitivity.

Finally, let \ts $\Des(w):= \{i\.:\. w(i)>w(i+1)\}$ \ts denote the \emph{set of descents} \ts in~$w$.
A permutation \ts $w\in S_n$ \ts is \emph{Grassmannian} \ts if it has at most one descent.

\smallskip

\subsection{Schubert coefficients}\label{ss:back-Sch}
Below we give a brief reminder of few basic results on Schubert polynomials.
We refer to \cite{Knu16,Mac91,Man01} for standard introductions to combinatorial
aspects, and to \cite{AF24,Ful97} for geometric aspects.

\defn{Schubert polynomials} \ts give a graded $\mathbb{Z}$-linear basis of polynomials ${\mathbb Z}[x_1,x_2,\ldots]$,
which we define recursively. Let $w_\circ=(n,n-1,\ldots,1) \in S_n\ts$.  Define
\begin{align*}
    {\mathfrak S}_{w_\circ}(x_1,\ldots,x_n)&:=x_1^{n-1}x_2^{n-2}\cdots x_{n-1}\,, \ \ \text{ and }\\
    {\mathfrak S_w}(x_1,\ldots,x_n)&:=\partial_i {\mathfrak S}_{ws_i}(x_1,\ldots,x_n) \ \ \text{ if } \ \ w(i)<w(i+1),
\end{align*}
where
$$\partial_i f \ := \ \frac{f-s_if}{x_i-x_{i+1}}\,.$$
Note that under the inclusion \.
$\iota: S_n\hookrightarrow S_{n+1}$ \. we have \.
${\mathfrak S}_w={\mathfrak S}_{\iota(w)}$.
This allows us to consider ${\mathfrak S}_w$ for each $w\in S_{\infty}$.

As Schubert polynomials \ts $\{{\mathfrak S}_w\}$ \ts form a polynomial basis, \defn{Schubert coefficients}
(also called \emph{Schubert structure coefficients}) \ts
$\{c_{u,v}^w\}$ \ts are defined by \eqref{eq:Schub-def}.
Although \ts $c_{u,v}^w\in\nn$,
they have no known combinatorial interpretation in full generality.
However, such an interpretation is known when \ts $v$ \ts is a simple transposition:

\smallskip

\begin{prop}[{\rm \defn{Monk's rule} \cite[$\S$2.7.1]{Man01}}{}] \label{prop:monk}
    Let \ts $u\in S_n$ \ts and let \ts $1\le k\le n-1$. Then:
    \[\mathfrak{S}_u \. \cdot \. \mathfrak{S}_{s_k}\ = \ \sum_{\substack{i\ts\leq \ts k\ts < \ts j\\
                 u \ts \lessdot \ts u t_{ij}}}\mathfrak{S}_{ut_{ij}}\..\]
\end{prop}

\smallskip

By their definition, Schubert polynomials \ts $\mathfrak{S}_{\rho}$ \ts
are homogeneous of degree \ts $\ell(\rho)$.  Thus we have:
\smallskip

\begin{prop}[{\rm \defn{dimension condition}}{}]\label{prop:dim}
    Suppose \. $\ell(w)\neq \ell(u)+\ell(v)$. \ts Then \. $c_{u,v}^w=0$.
\end{prop}
\smallskip

\smallskip

\subsection{Rothe diagrams}\label{sec:encode}
For a permutation \ts $w\in S_n\ts$, the \defn{Rothe diagram} \ts is defined as
\[
D(w) \, := \, \big\{(i,j): 1\leq i,j\leq n, \ts j<w(i), \ts i<w^{-1}(j)\big\} \, \subset \. \nn^2.
\]
Note that \ts $|D(w)|=\ell(w)$ \ts is the number of inversions in~$w$.
Depending on the context, we refer to elements of \ts $D(w)$ \ts as \emph{squares} \ts or \emph{boxes}.

For a permutation \ts $w\in S_\infty\ts$, its \emph{Lehmer code}  \ts
\ts $\code(w)=(c_1,c_2,\ldots)$ is such that \ts $c_i$ \ts is the number
of boxes in row~$i$ of $D(w)$.  We shorten this to ``code'' when the
context is clear.  For the identity permutation \ts ${\sf id}\in S_\infty\ts$, the
code is the all zero sequence.
Note that \ts $\code: S_\infty \to \kk$ \ts is a bijection, i.e.\
\ts $\code(w)$ \ts uniquely determines \ts $w\in S_{\infty}$ \ts and \ts $\code^{-1}$ \ts
is well defined.

\begin{ex}\label{ex:Rothe}
For a permutation
$w=72415836 \in {S}_{8}$, we have \ts $\code(w)=(6,1,2,0,1,2)$ \ts and \ts $\ell(w)=12$.
In this case \ts $n=8$, and the corresponding Rothe diagram \ts $D(w)$ \ts
is shown below. Here the squares are in positions \ts $(i,j)\in D(w)$, and the dots
are in positions \ts $(i,w(i))$, \ts $1\le i \le n$.

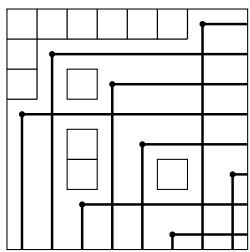
\[
\begin{tikzpicture}[scale=.4]
\draw (0,0) rectangle (8,8);

\draw[line width = .1ex] (0,7)--(6,7);
\draw[line width = .1ex] (1,8) -- (1,5);
\draw[line width = .1ex] (2,8) -- (2,7);
\draw[line width = .1ex] (3,8) -- (3,7);
\draw[line width = .1ex] (4,8) -- (4,7);
\draw[line width = .1ex] (5,8) -- (5,7);
\draw[line width = .1ex] (6,8) -- (6,7);

\draw[line width = .1ex] (0,6) -- (1,6);
\draw[line width = .1ex] (0,5) -- (1,5);

\draw (2,5) rectangle (3,6);

\draw (2,2) rectangle (3,4);
\draw[line width = .1ex] (2,3) -- (3,3);

\draw (6,3) rectangle (5,2);

\filldraw (6.5,7.5) circle (.5ex);
\draw[line width = .2ex] (6.5,0) -- (6.5,7.5) -- (8,7.5);
\filldraw (1.5,6.5) circle (.5ex);
\draw[line width = .2ex] (1.5,0) -- (1.5,6.5) -- (8,6.5);
\filldraw (3.5,5.5) circle (.5ex);
\draw[line width = .2ex] (3.5,0) -- (3.5,5.5) -- (8,5.5);
\filldraw (0.5,4.5) circle (.5ex);
\draw[line width = .2ex] (0.5,0) -- (0.5,4.5) -- (8,4.5);
\filldraw (4.5,3.5) circle (.5ex);
\draw[line width = .2ex] (4.5,0) -- (4.5,3.5) -- (8,3.5);
\filldraw (7.5,2.5) circle (.5ex);
\draw[line width = .2ex] (7.5,0) -- (7.5,2.5) -- (8,2.5);
\filldraw (2.5,1.5) circle (.5ex);
\draw[line width = .2ex] (2.5,0) -- (2.5,1.5) -- (8,1.5);
\filldraw (5.5,0.5) circle (.5ex);
\draw[line width = .2ex] (5.5,0) -- (5.5,0.5) -- (8,0.5);
\end{tikzpicture}
\]
\end{ex}

\medskip

\section{Code scaling}\label{s:Kir}

\subsection{Preliminaries}
\label{ss:Kir-prelim}
Let \ts $w\in S_n$ \ts be a permutation where \ts
$\code(w)=(c_1,c_2,\ldots,c_n)$.  Recall that for a given integer \ts
$N\ge 2$, the \defn{code scaling} \ts $N\ast w$ \ts is the
unique permutation with the code \ts $(N c_1,N c_2,\ldots,N c_n)$.

Observe that \ts $\ell(N\ast w)= N \ell(w)$.  Thus,
\ts $\ell(u)+\ell(v)=\ell(w)$ \ts if and only if \ts
$\ell(N\satA{}u)+\ell(N\satA{}v)=\ell(N\satA{}w)$ \ts for any $N>1$.
Additionally, it is clear by construction that code scaling preserves
the underlying descent sets of the permutations:  \ts $\Des(N\ast w) = \Des(w)$.

\begin{ex}\label{ex:Kir-prelim}
Take \ts $w= (7,2,4,1,5,8,3,6)$ \ts with \ts $\code(w)=(6,1,2,0,1,2)$ \ts
as in Example~\ref{ex:Rothe}.   Let \ts $N=2$.   By definition,
$$
2\ast w \, = \, \code^{-1}(12,2,4,0,2,4) \, =\, (13,3,6,1,5,9,2,4,7,8,10,11,12).
$$
We then have \ts $\ell(2\ast w) = 2\ts \ell(w)=24$ \ts and \ts $\Des(2\ast w) = \Des(w)=\{1,3,6\}$.
The corresponding Rothe diagram  $D(2\satA{} w)$ is the following:

 
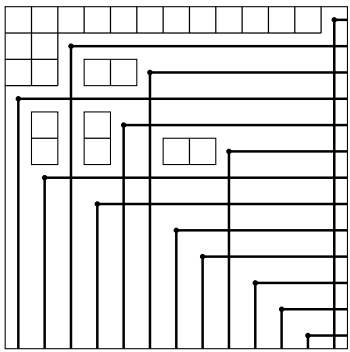
\[
\begin{tikzpicture}[scale=.35]
 \draw (0,0) rectangle (13,13);

\draw[line width = .1ex] (0,12)--(12,12);
\draw[line width = .1ex] (1,13) -- (1,10);
\draw[line width = .1ex] (2,10) -- (2,13);
\draw[line width = .1ex] (3,12) -- (3,13);
\draw[line width = .1ex] (4,12) -- (4,13);
\draw[line width = .1ex] (5,12) -- (5,13);
\draw[line width = .1ex] (6,12) -- (6,13);
\draw[line width = .1ex] (7,12) -- (7,13);
\draw[line width = .1ex] (8,12) -- (8,13);
\draw[line width = .1ex] (9,12) -- (9,13);
\draw[line width = .1ex] (10,12) -- (10,13);
\draw[line width = .1ex] (11,12) -- (11,13);
\draw[line width = .1ex] (12,12) -- (12,13);

\draw[line width = .1ex] (0,11)--(2,11);
\draw[line width = .1ex] (0,10)--(2,10);

\draw (3,10) rectangle (5,11);
\draw[line width = .1ex] (4,10) -- (4,11);

 \draw (1,7) rectangle (2,9);
\draw[line width = .1ex] (1,8) -- (2,8);

 \draw (3,7) rectangle (4,9);
\draw[line width = .1ex] (3,8) -- (4,8);

 \draw (6,7) rectangle (8,8);
\draw[line width = .1ex] (7,7) -- (7,8);

\filldraw (12.5,12.5) circle (.5ex);
\draw[line width = .2ex] (12.5,0) -- (12.5,12.5) -- (13,12.5);
\filldraw (2.5,11.5) circle (.5ex);
\draw[line width = .2ex] (2.5,0) -- (2.5,11.5) -- (13,11.5);
\filldraw (5.5,10.5) circle (.5ex);
\draw[line width = .2ex] (5.5,0) -- (5.5,10.5) -- (13,10.5);
\filldraw (0.5,9.5) circle (.5ex);
\draw[line width = .2ex] (0.5,0) -- (0.5,9.5) -- (13,9.5);
\filldraw (4.5,8.5) circle (.5ex);
\draw[line width = .2ex] (4.5,0) -- (4.5,8.5) -- (13,8.5);
\filldraw (8.5,7.5) circle (.5ex);
\draw[line width = .2ex] (8.5,0) -- (8.5,7.5) -- (13,7.5);
\filldraw (1.5,6.5) circle (.5ex);
\draw[line width = .2ex] (1.5,0) -- (1.5,6.5) -- (13,6.5);
\filldraw (3.5,5.5) circle (.5ex);
\draw[line width = .2ex] (3.5,0) -- (3.5,5.5) -- (13,5.5);
\filldraw (6.5,4.5) circle (.5ex);
\draw[line width = .2ex] (6.5,0) -- (6.5,4.5) -- (13,4.5);
\filldraw (7.5,3.5) circle (.5ex);
\draw[line width = .2ex] (7.5,0) -- (7.5,3.5) -- (13,3.5);
\filldraw (9.5,2.5) circle (.5ex);
\draw[line width = .2ex] (9.5,0) -- (9.5,2.5) -- (13,2.5);
\filldraw (10.5,1.5) circle (.5ex);
\draw[line width = .2ex] (10.5,0) -- (10.5,1.5) -- (13,1.5);
\filldraw (11.5,0.5) circle (.5ex);
\draw[line width = .2ex] (11.5,0) -- (11.5,0.5) -- (13,0.5);
\end{tikzpicture}
\]
\end{ex}

\smallskip

\subsection{St.~Dizier--Yong vanishing condition} \label{ss:Kir-SY}
For permutations \ts $u,v,w\in S_n\ts$, let \ts $\code(w)=(c_1,\ldots,c_n)$.
Consider integer fillings of boxes in \ts $D(u)\cup D(v)$ \ts
with entries in $[n]$. We view $D(u)\cup D(v)$ as a subset of \ts $[n]\times [2n]$,
where $D(v)$ is shifted right horizontally by $n$ units to the right of~$D(u)$.

In~\cite{ARY21}, Adve, the second author, and Yong define an
\defn{indicator tableau} \ts to be an integer filling \. $T:D(u)\cup D(v)\rightarrow [n]$,
such that
\begin{enumerate}
    \item the number of $i$'s in $T$ is equal to \ts $c_i= \code(w)_i$, for each \ts $i\in [n]$,
    \item each column of $T$ strictly increases from top to bottom, \ts and
    \item if an entry $m$ appears in row $r$ of $T$, then \ts $m\leq r$.
\end{enumerate}
Denote by \ts ${\sf Tab}_{u,v}^w$ \. to be the set of such indicator tableaux.
The \defn{St.~Dizier--Yong vanishing condition} \ts is the following:

\smallskip

\begin{lemma}[{\rm St.~Dizier--Yong \cite[Thm~B, $\S$4.3]{StDY22}}{}]\label{thm:stDY} \. We have:
$${\sf Tab}_{u,v}^w\. = \. \emp \ \ \Longrightarrow \ \ c_{u,v}^w \. = \. 0.
$$
\end{lemma}

\smallskip

\subsection{Proof of Theorem~\ref{t:sat-Kir}} \label{ss:Kir-proof}
For the first part, the equality \ts $c_{u,v}^w=1$ \ts
follows directly from Monk's rule (Proposition~\ref{prop:monk}) and
the construction of permutations $u,v,w$.

For the second part, let us show that the set of indicator tableaux
\ts ${\sf Tab}_{N\satA{}u,N\satA{}v}^{N\satA{}w}$ \ts is empty for each~$N>1$.
The result follows immediately by Lemma~\ref{thm:stDY}.

Let $S=\big\{c\in[u(i)+1,u(j)-1] \,: \, (j,c)\in D(u)\big\}$.
Since $u \lessdot u t_{ij}=w$, by the definition of the Rothe diagram,
    we have:
    \[D(w)\ = \ D(u) \. - \. \big\{(j,c) \,: \, c\in S\big\} \. \cup  \. \big\{(i,c) \,: \, c\in S\cup \{u(i)\}\big\}.\]
Let $\setsize=|S|$. By assumption (ii), $\al>0$. 
    Then we have
    \[\code(w) \, = \, \code(u)+(\setsize+1){\sf e}_i \. - \. (\setsize){\sf e}_j\.,
    \]
where \ts ${\sf e}_k = (0,\ldots,1,\ldots,0)$ \ts is the $k$-th standard basis vector.
    By definition of the code scaling, we then have:
    \begin{equation}\label{eq:code1}
        \code(N\satA{}w)\, = \, \code(N\satA{}u) \. + \. N(\setsize+1){\sf e}_i \. - \. N(\setsize){\sf e}_j\..
    \end{equation}

    Let \ts $D:=D(N\satA{}u)\cup D(N\satA{}v)$.  Suppose there exists an indicator tableau
\begin{equation}\label{eq:SY-tab}
    T\ts \in \ts {\sf Tab}_{N\satA{}u,N\satA{}v}^{N\satA{}w}
\end{equation}
    which by definition is a filling of~$D$.
    Note that for every $r\in[n]$, tableau $T$ contains \ts $ \code(N\satA{}w)_r$ \ts many $r$'s by assumption.
    Note also that every such entry of $r$ must appear in row $r$ or below.

\medskip

\nin
{\bf Claim:} \ts For all $(r,c)\in D$, we have:
    \[
T(r,c) \in \{i,j\} \ \, \text { if } \ \, r=j, \quad \text{and} \quad T(r,c) = r \ \ \text{otherwise}.
    \]

\begin{proof}
By \eqref{eq:code1}, we have \. $\code(N\satA{}w)_r=\code(N\satA{}u)_r$ \. for all \ts $r<i$.
This forces \ts $T(r,c)=r$ \ts for all \ts $(r,c)\in D$.
Similarly,  for \ts $r>j$, we have \. $\code(N\satA{}w)_r \ts=\ts\code(N\satA{}u)_r\ts$.
Thus again, we have \ts $T(r,c)=r$ \ts for all \ts $(r,c)\in D$.

For the case of $r=i$, Equation~\eqref{eq:code1} gives:
\[
\code(N\satA{}w)_i \, = \, \code(N\satA{}u)_i \. + \. N(\setsize+1).
\]
Similarly, this implies that \ts $T(i,c)=i$ \ts for all \ts $(i,c)\in D$.
Note that \ts $D(N\satA{}v)$ \ts contains $N$ boxes in row $i$ and no boxes elsewhere.
This leaves
\[
\big(N\code(u)_i+N(\setsize+1)\big) \. - \. \big(N\code(u)_i+N\code(v)_i\big) \, = \, N \setsize \, > \, 0
\]
many $i$'s left to place in $T$.

For the case \ts $i<r<j$, we again have \.
$ \code(N\satA{}w)_r=\code(N\satA{}u)_r\ts$. Since
$i<r<j$, by assumptions (i) and (iii) we have $u(r)<u(i)$.
Note that any permutation $u'\in S_n$ and $1\leq i<j\leq n$, the following are equivalent by the definition of Rothe diagrams:
\begin{enumerate}
    \item[(a)] For all $i<r<j$, $u'(r)<u'(i)$, and
    \item[(b)] For all $i<r<j$, $\code(u')_r\leq \code(u')_i-(r-i)$.
\end{enumerate}
Additionally, (a) or (b) implies that 
\begin{enumerate}
    \item[(c)] For all $i<r<j$, if $(r,c)\in D(u')$, then $(i,c)\in D(u')$.
\end{enumerate}

Since by assumption, $u$ satisfies (a), we know $u$ satisfies (b). Since $\code(N\satA{}u)=N\code(u)$, we know $N\satA{}u$ satisfies (b), so $N\satA{}u$ satisfies (a).  Therefore, both $u$ and $N\satA{}u$ have property (c).

Since there can be no repetition
of entries in columns of~$T$, this gives $T(r,c)\neq i$ for $i<r<j$. Then similarly as in the previous cases, this forces
\ts $T(r,c)=r$ \ts for $(r,c)\in D$.
Finally, for the case $r=j$, this leaves \ts $T(j,c)\in\{i,j\}$ \ts for all \ts $(j,c)\in D$.
\end{proof}

\smallskip

Now we focus on the entries in row $j$ of $T$.
By the Claim, this row must contain the remaining \ts $N \setsize$ \ts entries~$i$.
Since the entries in columns of $T$ must be distinct, and row $i$ contains all $i$'s,
     these \ts $N \setsize$ \ts many
     $i$'s cannot appear below a square from row~$i$.

    Recall that row $j$ of $D$ contains
     $\code(N\satA{}u)_j=N\code(u)_j$ many squares.
    If we have that $N\code(u)_i>N\code(u)_j$, then (a) holds, so (c) holds. Thus all boxes in row $j$ would lie beneath boxes in row $i$, leaving no available boxes for these $N\setsize$ remaining $i$'s, a contradiction. 

Thus we assume $N\code(u)_i<N\code(u)_j$.
    Then similarly as above,
    it follows that $D$
    has
    \[
    \code(N\satA{}u)_i-(j-i-1) \, = \, N\code(u)_i-(j-i-1)
    \]
    many squares in row $j$ of $D$ directly below boxes in row $i$ of $D$.
    Summarizing, since $T$ must have distinct column entries, we need:
    \begin{equation}\label{eq:codeReq}
        N\code(u)_j \. - \. ( N\code(u)_i-(j-i-1)) \, \geq \, N\setsize.
    \end{equation}

     By~\eqref{eq:code1} and the fact that \ts $u(k)<u(i)$ for each $i<k<r$, we have:
    \begin{equation*}
        \code(u)_j \, = \, (\code(u)_i-(j-i-1)) \. + \. \setsize.
    \end{equation*}
    Applying this to the left-hand side of Equation~\eqref{eq:codeReq} gives
    \begin{equation}\label{eq:code2N}
        N\code(u)_j \. - \. \big(N\code(u)_i-(j-i-1)\big) \, = \, N\setsize\. - \. (N-1)(j-i-1).
    \end{equation}
    Combining Equations~\eqref{eq:codeReq} and~\eqref{eq:code2N}, we conclude:
    \[0\, \geq \, (N-1)(j-i-1),
    \]
    a contradiction with assumptions \ts $j-i\geq 2$ \ts and \ts $N> 1$.
We conclude that there is no $T$ as in \eqref{eq:SY-tab}, i.e.\
\. ${\sf Tab}_{N\satA{}u,N\satA{}v}^{N\satA{}w}=\emp$.
The result now follows by Lemma~\ref{thm:stDY}. \qed

\smallskip

\subsection{Examples and special cases} \label{ss:Kir-ex}
We start with the simplest example which both illustrates the proof above
and motivates other special cases.

\smallskip

\begin{ex}\label{ex:star}
   Take \ts $u=2143$, \ts $v=s_1=2134$ \ts and \ts $w=ut_{13}=4123$.
   Note that \. $\code({u})=(1,0,1)$, \. $\code({v})=(1)$, \.
   $\code({w})=(3)$, and that \ts $u\lessdot w$.
By Monk's rule, we have \ts $c_{u,v}^{w}=1$.

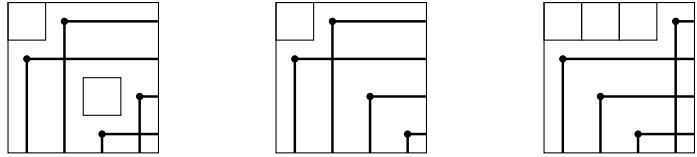
\[
\begin{tikzpicture}[scale=0.5]
\draw (0,0) rectangle (4,4);

 \draw (0,4) rectangle (1,3);
\draw (2,2) rectangle (3,1);

\filldraw (1.5,3.5) circle (.5ex);
\draw[line width = .2ex] (1.5,0) -- (1.5,3.5) -- (4,3.5);
\filldraw (0.5,2.5) circle (.5ex);
\draw[line width = .2ex] (0.5,0) -- (0.5,2.5) -- (4,2.5);
\filldraw (3.5,1.5) circle (.5ex);
\draw[line width = .2ex] (3.5,0) -- (3.5,1.5) -- (4,1.5);
\filldraw (2.5,0.5) circle (.5ex);
\draw[line width = .2ex] (2.5,0) -- (2.5,0.5) -- (4,0.5);

\end{tikzpicture}
\qquad \qquad
\begin{tikzpicture}[scale=0.5]
\draw (0,0) rectangle (4,4);

 \draw (0,4) rectangle (1,3);

\filldraw (1.5,3.5) circle (.5ex);
\draw[line width = .2ex] (1.5,0) -- (1.5,3.5) -- (4,3.5);
\filldraw (0.5,2.5) circle (.5ex);
\draw[line width = .2ex] (0.5,0) -- (0.5,2.5) -- (4,2.5);
\filldraw (2.5,1.5) circle (.5ex);
\draw[line width = .2ex] (2.5,0) -- (2.5,1.5) -- (4,1.5);
\filldraw (3.5,0.5) circle (.5ex);
\draw[line width = .2ex] (3.5,0) -- (3.5,0.5) -- (4,0.5);
\end{tikzpicture}
\qquad\qquad
\begin{tikzpicture}[scale=0.5]
\draw (0,0) rectangle (4,4);

 \draw (0,4) rectangle (1,3);
  \draw (2,4) rectangle (1,3);
    \draw (2,4) rectangle (3,3);

\filldraw (3.5,3.5) circle (.5ex);
\draw[line width = .2ex] (3.5,0) -- (3.5,3.5) -- (4,3.5);
\filldraw (0.5,2.5) circle (.5ex);
\draw[line width = .2ex] (0.5,0) -- (0.5,2.5) -- (4,2.5);
\filldraw (1.5,1.5) circle (.5ex);
\draw[line width = .2ex] (1.5,0) -- (1.5,1.5) -- (4,1.5);
\filldraw (2.5,0.5) circle (.5ex);
\draw[line width = .2ex] (2.5,0) -- (2.5,0.5) -- (4,0.5);

\end{tikzpicture}
\]

Observe that \. $2\satA{}u=31524$, \. $2\satA{}v=21345=v$, and \. $2\satA{}w=7123456$.
Since \. $\code({2\satA{}w})=(6)$, there is a unique way to fill \. $D=D(2\satA{}u) \cup D(2\satA{}v)$ \.
with $1$'s.

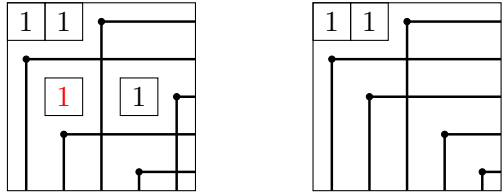
\[
\begin{tikzpicture}[scale=0.5]
\draw (0,0) rectangle (5,5);

\draw (0,4) rectangle (1,5) node[pos=.5] {$1$};
\draw (1,4) rectangle (2,5) node[pos=.5] {$1$};

\draw (1,2) rectangle (2,3) node[pos=.5] {\textcolor{red}{$1$}};

\draw (3,2) rectangle (4,3) node[pos=.5] {$1$};

\filldraw (2.5,4.5) circle (.5ex);
\draw[line width = .2ex] (2.5,0) -- (2.5,4.5) -- (5,4.5);
\filldraw (0.5,3.5) circle (.5ex);
\draw[line width = .2ex] (0.5,0) -- (0.5,3.5) -- (5,3.5);
\filldraw (4.5,2.5) circle (.5ex);
\draw[line width = .2ex] (4.5,0) -- (4.5,2.5) -- (5,2.5);
\filldraw (1.5,1.5) circle (.5ex);
\draw[line width = .2ex] (1.5,0) -- (1.5,1.5) -- (5,1.5);
\filldraw (3.5,0.5) circle (.5ex);
\draw[line width = .2ex] (3.5,0) -- (3.5,0.5) -- (5,0.5);

\end{tikzpicture}
\qquad\qquad
\begin{tikzpicture}[scale=0.5]
\draw (0,0) rectangle (5,5);

\draw (0,4) rectangle (1,5) node[pos=.5] {$1$};
\draw (1,4) rectangle (2,5) node[pos=.5] {$1$};

\filldraw (2.5,4.5) circle (.5ex);
\draw[line width = .2ex] (2.5,0) -- (2.5,4.5) -- (5,4.5);
\filldraw (0.5,3.5) circle (.5ex);
\draw[line width = .2ex] (0.5,0) -- (0.5,3.5) -- (5,3.5);
\filldraw (1.5,2.5) circle (.5ex);
\draw[line width = .2ex] (1.5,0) -- (1.5,2.5) -- (5,2.5);
\filldraw (3.5,1.5) circle (.5ex);
\draw[line width = .2ex] (3.5,0) -- (3.5,1.5) -- (5,1.5);
\filldraw (4.5,0.5) circle (.5ex);
\draw[line width = .2ex] (4.5,0) -- (4.5,0.5) -- (5,0.5);

\end{tikzpicture}
\]

By the definition of indicator tableaux, they must be increasing in columns,
ruling out the filling above.  Thus \. ${\sf Tab}_{2\satA{}u,2\satA{}v}^{2\satA{}w} =\emp$.
By Lemma~\ref{thm:stDY}, we have \. $c_{2\satA{}u,2\satA{}v}^{2\satA{}w}=0$,
giving a counterexample to the saturation property \eqref{eq:sat-Kir}.
\end{ex}

\smallskip

The following corollary follows immediately from Theorem~\ref{t:sat-Kir}
by taking \ts $j=i+2$.

\smallskip

\begin{cor}\label{c:sat-Kir-gen}
Let \. $u\in S_n$ \. such that \. $u(i+1)<u(i)< u(i+2)-1$, and additionally that $\{c \. : \. u(i)<c<u(i+2), \. i+2<u^{-1}(c) \}\neq\emp$.
Let \. $v=(i,i+1)$ \. and let \. $w=u \cdot (i,i+2)$. Then:
\begin{equation}\label{eq:sat-Kir-gen}
c_{u,v}^w\.= \. 1 \qquad \text{and}  \qquad c_{N\satA{}u,N\satA{}v}^{N\satA{}w}\. = \. 0 \ \ \text{for all}  \ \ N>1.
\end{equation}
\end{cor}
\smallskip

Note that Corollary~\ref{c:main-Kir} is a special case
of Corollary~\ref{c:sat-Kir-gen}.  Note also that
Example~\ref{ex:star} is a special case of Corollary~\ref{c:main-Kir}
when \ts $n=4$.





\medskip
\section{Final remarks}\label{s:finrem}

\subsection{}\label{ss:finrem-CS}
Our own motivation to study the saturation property for Schubert
coefficients lies in connection to the \defn{Schubert vanishing problem}
\ts $\big[\ts c^w_{u,v} \ts =^? \ts 0 \ts\big]$.  The idea here is
to extend the approach in \cite{DM06,MNS12} to the
\defn{LR vanishing problem} \ts
$\big[\ts c^\la_{\mu\nu} \ts =^? \ts 0 \ts\big]$.
There, the authors independently observed\footnote{The original preprints appeared on the
{\tt arXiv} in January 2005, within a day of each other.} that the
saturation property \eqref{eq:saturation} implies that the vanishing
of LR--coefficients can be solved by a linear program (LP).  This gives
a deterministic poly-time algorithm for deciding the LR~vanishing.

There are two main ingredients in the approach above, both representing
a major obstacle.
First, one needs a combinatorial interpretation of the LR~coefficient
\ts $c^\la_{\mu\nu}$ \ts to show that it counts the number of integer
points in a convex polytope \ts $Q_{\la\mu\nu}$ \ts defined by integer constraints.
Second, one needs a saturation property to reduce the LR~vanishing problem
to \ts $Q_{\la\mu\nu}$ \ts containing a \emph{rational} \ts point.  One then
applies known results that LP is in~$\poly$ to conclude the same for the
LR~vanishing.

In our most recent paper \cite{PR-BPP} capping a series of weaker
results, we show that the vanishing of Schubert coefficients \ts
$\big[\ts c^u_{vw} \ts =^? \ts 0 \ts\big]$ \ts is in \ts $\coRP\subseteq \BPP$,
i.e.\ can be decided in probabilistic polynomial time with a one-sided error
(in the case of a positive answer).  This is very low in the polynomial hierarchy,
and suggests a possibility that Schubert vanishing might be in~$\poly$,
at least for some classes of permutations.

Since the reduction to LP is really the only approach that we know
(short of derandomization), we would need the two ingredients
described above.  As mentioned in the introduction,  there
are at least two different combinatorial interpretations of
\emph{$2$-step Schubert coefficients}, which correspond to
permutations with at most two descents \cite{Cos09,BKPT16}.
While we are not aware if there is a way to restate either of these combinatorial
interpretations in terms of the number of integer points in polytopes,
this case is a natural place to start.

\begin{question}\label{q:2-step}
Can the vanishing problem for $2$-step Schubert coefficients be decided
in poly-time?
\end{question}

Unfortunately, our Corollary~\ref{c:main-Kir}
implies that in the $2$-step case, the natural saturation property fails.
The question remains wide open, unlikely to be resolved by the existing tools.

\subsection{} \label{ss:finrem-sat-other}
In connection to Question~\ref{q:2-step}, it is worth comparing how
the failure of saturation properties for other algebraic combinatorics
constants discussed in the introduction
relate to the complexity of the corresponding vanishing problems.

For example, it remains open whether the vanishing of the \emph{Clebsch--Gordan} (CG)
\emph{coefficients} (particular generalizations of LR coefficients to other root systems)
can be decided in poly-time.  The problem is resolved for the even weights
by a combination of results in \cite{DM06} and \cite{KM08}.  It is known that
the Ehrhart positivity conjectures
in \cite{KTT04} would imply that  \emph{CG~vanishing} is in~$\poly$ in full
generality.  Unfortunately, these conjectures remain wide open even in the
special case of Kostka numbers, see e.g.\ \cite{Ale19}.

In some cases, there are known computational complexity obstacles for the vanishing
problem, immediately invalidating a saturation property approach as above.
Notably, it is known that the vanishing problem for the Kronecker
coefficients is $\NP$-hard \cite{IMW17}.  Later, in \cite{PP20}, it was shown
that the vanishing of reduced Kronecker coefficients is also $\NP$-hard.
By itself, the $\NP$-hardness of the vanishing problem does not automatically imply
that the saturation property fails, but it does suggest a more involved
underlying structure in the problem.

\subsection{}\label{ss:finrem-other}
Finally, note that there is more than one way to define an operation on permutations
to produce a saturation property.  In the previous {\tt arXiv} version ({\tt v1}) of this paper \cite[$\S$4]{PR-Sat},
we obtain similar results for the folklore \emph{bit scaling} \ts operation,
which can be viewed as a partial tensor product with the identity permutation.
The saturation property in this case also generalizes the saturation property
for LR~coefficients, so it is natural to ask whether this property holds.
We refute this possibility in a strong sense; the proof in this case uses Monk's rule
and the \emph{dimension condition}.

There are other possible scaling operations which can be defined on permutations,
i.e.\ by replacing each box of the Rothe diagram with an $N \times N$ square of boxes.
Unfortunately, code scaling and bit scaling are the only operations we know
that preserve the descents. This is a necessary condition to ensure that
the operation applied to Grassmannian permutations corresponds to the usual
partition multiplication.

\vskip.4cm

\subsection*{Acknowledgements}
We thank Sara Billey, Jes\'us De Loera, Pavel Galashin, Oliver Pechenik,
Zach Slonim, Frank Sottile, Alex Yong, Sylvester Zhang and Paul Zinn-Justin
for interesting discussions and helpful comments. We also thank the referees for their careful feedback and suggested improvements.


\vskip.7cm

\newpage

{\footnotesize

\begin{thebibliography}{abcdefgh}


\bibitem[ARY19]{ARY19}
Anshul~Adve, Colleen~Robichaux and Alexander~Yong,
Vanishing of Littlewood--Richardson polynomials is in $\poly$,
\emph{Comput.\ Complexity}~\textbf{28} (2019), 241--257.

\bibitem[ARY21]{ARY21}
Anshul~Adve, Colleen~Robichaux and Alexander~Yong,
An efficient algorithm for deciding vanishing of Schubert polynomial coefficients,
\emph{Adv.\ Math.}~\textbf{383} (2021), Paper No.~107669, 38~pp.;
extended abstract in \emph{Proc.\ 31st FPSAC} (2020), Art.~52, 12~pp.




\bibitem[Ale19]{Ale19}
Per~Alexandersson,
Polytopes and large counterexamples,
\emph{Exp.\ Math.}~\textbf{28} (2019), 115--120.

\bibitem[AF24]{AF24}
David~Anderson and William~Fulton,
\emph{Equivariant cohomology in algebraic geometry},
Cambridge Univ.\ Press, Cambridge, UK, 2024, 446~pp.

\bibitem[ARY13]{ARY13}
David~Anderson, Edward~Richmond and Alexander~Yong,
Eigenvalues of Hermitian matrices and equivariant cohomology of Grassmannians,
\emph{Compos.\ Math.}~\textbf{149} (2013), 1569--1582.





\bibitem[Bel06]{Bel06}
Prakash~Belkale, Geometric proofs of {H}orn and saturation conjectures,
\emph{J.\ Algebraic Geom.}~\textbf{15} (2006), 133--173.

\bibitem[Bel08]{Bel08}
Prakash~Belkale, Quantum generalization of the Horn conjecture,
\emph{Jour.\  AMS}~\textbf{21} (2008), 365--408.

\bibitem[BK10]{BK10}
Prakash~Belkale and Shrawan~Kumar,
Eigencone, saturation and Horn problems for symplectic and odd orthogonal groups,
\emph{J.~Algebraic Geom.} \textbf{19} (2010), 199--242.

\bibitem[BKR12]{BKR12}
Prakash~Belkale, Shrawan~Kumar and Nicolas~Ressayre,
A generalization of Fulton's conjecture for arbitrary groups,
\emph{Math.\ Ann.}~\textbf{354} (2012), 401--425.




\bibitem[BVW17]{BVW17}
Nicole~Berline, Mich\`ele Vergne and Michael Walter,
The Horn inequalities from a geometric point of view, \emph{Enseign.\ Math.}~\textbf{63} (2017), 403--470.




\bibitem[BOR09]{BOR09}
Emmanuel~Briand, Rosa~Orellana and Mercedes~Rosas,
Reduced Kronecker coefficients and counter-examples to Mulmuley's
strong saturation conjecture SH (with an appendix by K.~Mulmuley),
\emph{Comp.\ Complexity}~\textbf{18} (2009), 577--600.


\bibitem[Buch00]{Buch00}
Anders~S.~Buch,
The saturation conjecture (after A.~Knutson and T.~Tao).
With an appendix by William~Fulton,
\emph{Enseign.\ Math.}~\textbf{46} (2000), 43--60.

\bibitem[Buch02]{Buch02}
Anders~S.~Buch,
A Littlewood--Richardson rule for the $K$-theory of Grassmannians,
\emph{Acta Math.}~\textbf{189} (2002), 37--78.

\bibitem[BKPT16]{BKPT16}
Anders~S.~Buch,   Andrew~Kresch, Kevin~Purbhoo and Harry~Tamvakis,
The puzzle conjecture for the cohomology of two-step flag manifolds,
\emph{J.~Algebraic Combin.} \textbf{44} (2016), 973--1007.

\bibitem[CR23]{CR23}
Pierre-Emmanuel~Chaput and Nicolas~Ressayre,
Reduction for branching multiplicities,
\emph{Int.\ Math.\ Res.\ Not.\ IMRN} \textbf{2023} (2023), 15207--15265.

\bibitem[CHM07]{CHM07}
Matthias~Christandl, Aram~W.~Harrow and Graeme~Mitchison, 
Nonzero Kronecker coefficients and what they tell us about spectra,
\emph{Comm.\ Math.\ Phys.}~\textbf{270} (2007), 575--585.

\bibitem[Cos09]{Cos09}
Izzet~Coskun,
A Littlewood--Richardson rule for two-step flag varieties,
\emph{Invent.\ Math.}~\textbf{176} (2009), 325--395.

\bibitem[DM06]{DM06}
Jes\'us~A.~De~Loera and Tyrrell~B.~McAllister, On the computation of
{C}lebsch--{G}ordan coefficients and the dilation effect,
\emph{Experiment.\ Math.}~\textbf{15} (2006), 7--19.

\bibitem[DW00]{DW00}
Harm~Derksen and Jerzy~Weyman,
Semi-invariants of quivers and saturation for Littlewood--Richardson coefficients,
\emph{Jour.~AMS} \textbf{13} (2000), 467--479.

\bibitem[\`Ela92]{Ela92}
Alexander~G.~\`Elashvili,
Invariant algebras, in \emph{Lie groups, their discrete subgroups,
and invariant theory}, AMS, Providence, RI, 1992, 57--64.

\bibitem[Ful97]{Ful97}
William~Fulton, \emph{Young tableaux},
Cambridge Univ.\ Press, Cambridge, UK, 1997, 260~pp.

\bibitem[Ful00]{Ful00}
William~Fulton,
Eigenvalues, invariant factors, highest weights, and Schubert calculus,
\emph{Bull.\ AMS}~\textbf{37} (2000), 209--249.

\bibitem[GOY21]{GOY21}
Shiliang~Gao, Gidon~Orelowitz and Alexander~Yong,
Newell--Littlewood numbers, \emph{Trans.~AMS} \textbf{374} (2021), 6331--6366.

\bibitem[IMW17]{IMW17}
Christian~Ikenmeyer, Ketan~D.~Mulmuley and Michael~Walter,
On vanishing of Kronecker coefficients,
\emph{Comp.\ Complexity}~\textbf{26} (2017), 949--992.

\bibitem[KM08]{KM08}
Michael~Kapovich and John~J.~Millson,
A path model for geodesics in Euclidean buildings and its applications
to representation theory, \emph{Groups Geom.\ Dyn.}~\textbf{2}
(2008), 405--480.

\bibitem[KTT04]{KTT04}
Ronald~C.~King, Christophe~Tollu and Fr\'ed\'eric~Toumazet,
Stretched Littlewood--Richardson and Kostka coefficients,
in \emph{Symmetry in physics}, AMS, Providence, RI, 2004, 99--112.

\bibitem[Kir04]{Kir04}
Anatol~N.~Kirillov,
An invitation to the generalized saturation conjecture,
\emph{Publ.\ RIMS}~\textbf{40} (2004), 1147--1239.

\bibitem[Kly04]{Kly04}
Alexander~Klyachko,
Quantum marginal problem and representations of the symmetric group,
preprint (2004), 47~pp.;
{\tt arXiv:quant-ph/0409113}.

\bibitem[Knu16]{Knu16}
Allen~Knutson, Schubert calculus and puzzles, in
\emph{Adv.\ Stud.\ Pure Math.}~\textbf{71}, Math.\ Soc.\ Japan, Tokyo, 2016, 185--209.

\bibitem[Knu22]{Knu22}
Allen~Knutson, Schubert calculus and quiver varieties, in
\emph{Proc.\ ICM} (2022, virtual), Vol.~VI, EMS Press, 4582--4605.

\bibitem[KT99]{KT99}
Allen~Knutson and Terence~Tao,
The honeycomb model of $\GL_n(\cc)$ tensor products~I:
Proof of the saturation conjecture,
{\em Jour.~AMS}~\textbf{12} (1999), 1055--1090.

\bibitem[KT03]{KT03}
Allen~Knutson and Terence~Tao, Puzzles and (equivariant) cohomology of
Grassmannians, \emph{Duke Math.~J.} \textbf{119} (2003), 221--260.

\bibitem[KTW04]{KTW04}
Allen~Knutson, Terence~Tao, and Christopher~Woodward. The honeycomb model of $GL_n(\mathbb{C})$ tensor products II: Puzzles determine facets of the Littlewood-Richardson cone. \emph{Jour.\  AMS}, 17(1), (2004) 19--48.

\bibitem[Kum14]{Kum14}
Shrawan~Kumar,
A survey of the additive eigenvalue problem.
With an appendix by M.~Kapovich,
\emph{Transform.\ Groups} \textbf{19} (2014), 1051--1148.

\bibitem[Mac91]{Mac91}
Ian~G.~Macdonald,
\emph{Notes on Schubert polynomials},
Publ.\ LaCIM, UQAM, Montreal, 1991, 116~pp.; available at \ts \href{http://www.math.uwaterloo.ca/~opecheni/macdonaldschubert.pdf}{tinyurl.com/382f7an7}





\bibitem[Mac95]{Mac95}
Ian~G.~Macdonald, \emph{Symmetric functions and Hall polynomials}
(Second ed.), Oxford U.~Press, New York, 1995, 475~pp.

\bibitem[Man01]{Man01}
Laurent~Manivel,
\emph{Symmetric functions, Schubert polynomials and degeneracy loci},
SMF/AMS, Providence, RI, 2001, 167~pp.



\bibitem[Min24]{Min24}
Jaewon~Min,
Proof of the Newell--Littlewood saturation conjecture, preprint (2024), 50 pp.;
{\tt arXiv:2409.} {\tt 00233}.

\bibitem[Mul09]{Mul09}
Ketan~D.~Mulmuley,
Geometric Complexity Theory~VI: the flip via saturated and positive integer programming
in representation theory and algebraic geometry, preprint (2009, v4), 139~pp.; \ts
{\tt arXiv:} {\tt 0704.0229}.

\bibitem[MNS12]{MNS12}
Ketan~D.~Mulmuley, Hariharan~Narayanan and Milind~Sohoni,
Geometric complexity theory~III.  On deciding nonvanishing of a
Littlewood--Richardson coefficient, \emph{J.\ Algebraic Combin.}~\textbf{36}
(2012), 103--110.




\bibitem[PP20]{PP20}
Igor~Pak and Greta~Panova,
Breaking down the reduced Kronecker coefficients, \emph{C.~R.\ Math.\ Acad.\ Sci.\ Paris}
\textbf{358} (2020), no.~4, 463--468.



\bibitem[PR25]{PR-BPP}
Igor~Pak and Colleen~Robichaux,
Vanishing of Schubert coefficients in probabilistic polynomial time,
preprint (2025), 15~pp.; \ts {\tt arXiv:2509.16467}.


\bibitem[PR26]{PR-Sat}
Igor~Pak and Colleen~Robichaux,
Saturation property fails for Schubert coefficients, preprint (2026), 13~pp.; \ts {\tt arXiv:2601.04182}, {\tt v1}.


\bibitem[PS26]{PS26}
Igor~Pak and Zachary~Slonim, 
Stretched Schubert coefficients are eventually quasi-polynomial, 
preprint (2026), 25~pp.; \ts {\tt arXiv:2604.27107}.



\bibitem[RYY22]{RYY22}
Colleen~Robichaux, Harshit~Yadav and Alexander Yong,
Equivariant cohomology, Schubert calculus, and edge labeled tableaux,
in \emph{Facets of Algebraic Geometry}, Vol.~II,
Cambridge Univ.\ Press, Cambridge, UK, 2022, 284--335.

\bibitem[Ros01]{Ros01}
Mercedes~H.~Rosas,  The Kronecker product of Schur
functions indexed by two-row shapes or hook shapes,
\emph{J.\ Algebraic Combin.} \textbf{14} (2001), 153--173.




\bibitem[SS16]{SS16}
Steven~V.~Sam and Andrew~Snowden,
Proof of Stembridge's conjecture on stability of Kronecker coefficients,
\emph{J.~Algebraic Combin.}~\textbf{43} (2016), 1--10.

\bibitem[SY22]{StDY22}
Avery~St.~Dizier and Alexander~Yong, Generalized permutahedra and
{S}chubert calculus, \emph{Arnold Math.~J.} \textbf{8} (2022), 517--533.

\bibitem[Sta99]{Sta99}
Richard~P.~Stanley, {\em Enumerative Combinatorics}, vol.~1 (Second ed.)
and vol.~2, Cambridge Univ.~Press, 2012 and~1999, 626~pp.\ and 581~pp.




\bibitem[Zel99]{Zel99}
Andrei~Zelevinsky,
Littlewood--Richardson semigroups, in \emph{New perspectives in algebraic combinatorics},
Cambridge Univ.\ Press, Cambridge, UK, 1999, 337--345.

\end{thebibliography}
}
\end{document}